\documentclass[final,onefignum,onetabnum]{siamart220329}

\usepackage{lipsum}
\usepackage{amsfonts}
\usepackage{graphicx}
\usepackage{epstopdf}
\usepackage{algorithmic}
\ifpdf
  \DeclareGraphicsExtensions{.eps,.pdf,.png,.jpg}
\else
  \DeclareGraphicsExtensions{.eps}
\fi

\usepackage[clean]{revdiff}

\headers{Falling liquid film control}{O. A. Holroyd, R. Cimpeanu, and S. N. Gomes}

\title{Linear quadratic regulation control for falling liquid films\thanks{Submitted to the editors January 26 2023.
\funding{This work was funded by the UK Engineering and Physical Sciences Research Council (EPSRC) through grants EP/S022848/1 and EP/V051385/1.}}}

\author{Oscar A. Holroyd\thanks{Mathematics Institute, University of Warwick, Coventry CV4 7AL, UK
  (\email{o.holroyd@warwick.ac.uk}, \email{radu.cimpeanu@warwick.ac.uk}, \email{susana.gomes@warwick.ac.uk}).}
\and Radu Cimpeanu\footnotemark[2]
\and Susana N. Gomes\footnotemark[2]}

\ifpdf
\hypersetup{
  pdftitle={Falling Film LQR Control},
  pdfauthor={O. A. Holroyd, R. Cimpeanu, and S. N. Gomes}
}
\fi

\usepackage{todonotes}
\makeatletter
\@mparswitchfalse%
\makeatother
\normalmarginpar 

\usepackage{amsopn}
\DeclareMathOperator{\tr}{tr} 

\newcommand*{\tran}{^{\mkern-1.5mu\mathsf{T}}} 
\newcommand{\abs}[1]{\left\lvert#1\right\rvert} 
\renewcommand{\vec}[1]{\textbf{#1}} 
\renewcommand{\matrix}[1]{ 
  \begin{pmatrix}
    #1
  \end{pmatrix}
}
\newcommand*{\unit}[1]{\text{\hspace{0.15em}#1}}
\newcommand*{\diff}{\mathop{}\!\mathrm{d}}

\newcommand{\Rey}{\mbox{\textit{Re}}}
\newcommand{\Ca}{\mbox{\textit{Ca}}}

\begin{document}
\maketitle

\begin{abstract}
We propose and analyse a new methodology based on linear-quadratic regulation (LQR) for stabilising falling liquid films via blowing and suction at the base. LQR methods enable rapidly responding feedback control by precomputing a gain matrix, but are only suitable for systems of linear ordinary differential equations (ODEs). By contrast, the Navier-Stokes equations that describe the dynamics of a thin liquid film flowing down an inclined plane are too complex to stabilise with standard control-theoretical techniques. To bridge this gap we use reduced-order models -- the Benney equation and a weighted-residual integral boundary layer model -- obtained via asymptotic analysis to derive a multi-level control framework. This framework consists of an LQR feedback control designed for a linearised and discretised system of ODEs approximating the reduced-order system, which is then applied to the full Navier-Stokes system. The control scheme is tested via direct numerical simulation (DNS), and compared to analytical predictions of linear stability thresholds and minimum required actuator numbers. Comparing the strategy between the two reduced-order models we show that in both cases we can successfully stabilise towards a uniform flat film across their respective ranges of valid parameters, with the more accurate weighted-residual model outperforming the Benney-derived controls. The weighted-residual controls are also found to work successfully far beyond their anticipated range of applicability. The proposed methodology increases the feasibility of transferring robust control techniques towards real-world systems, and is also generalisable to other forms of actuation.
\end{abstract}



\section{Introduction}
  \label{sec:introduction}
  Modelling and stabilisation of falling liquid films is a fundamental problem at the intersection of fluid dynamics, asymptotic analysis, and control theory. Manipulation of these multi-scale systems is key to a number of industrial applications ranging from coating flows in liquid crystal display devices to microchip manufacture. Such systems have a high degree of complexity, which makes them challenging to model, control and simulate accurately and efficiently.

  Although the control of complex systems is a challenge found in a wide range of applied sciences -- from preventing ice buildup on aerofoils~\cite{palacios2011instantaneous} to avoiding obstacles in self-driving vehicles~\cite{paden2016survey} and crowd management~\cite{burger2013mean} -- falling liquid films are a prototypical example of such a control problem. Governed by the two-phase Navier-Stokes equations, these are multi-scale setups used in many applications as well as beautiful day-to-day phenomena such as wavy films on a window on a rainy day. The resulting flow becomes unstable above a critical Reynolds number (a parameter depending on velocity, inclination angle, film thickness, and fluid density), exhibiting a rich set of behaviours starting with two-dimensional (2D) waves and leading to 3D spatiotemporal chaos. \rnew{This has been shown experimentally by Nosoko et al.~\cite{nosoko996characteristics}, and analytically for simplified models by Chat\'e and Manneville~\cite{chate1987transition} and Smyrlis and Papageorgiou~\cite{smyrlis1991predicting}.} Since the 1960s the problem has attracted much analytical focus, resulting in a number of periodic reduced-order models~\cite{benney1966long,ruyer2000improved,sivashinsky1980irregular,wray2020reduced} based on the assumption that perturbations are `long-wave', i.e.\ their wavelength is much larger than their amplitude. Such models range from single equation models describing the dynamics of the liquid film height, to systems of multiple equations describing the height, downstream flux, and potentially additional independent quantities. A broad range of thin-film models are covered extensively in two reviews by Craster and Matar~\cite{craster2009dynamics} and Kalliadasis et al.~\cite{kalliadasis2011falling}, and these models were recently coupled with validation in experimental settings by Denner et al.~\cite{denner2018solitary}. More recently, Richard et al~\cite{richard2019optimization}, Usha, Chattopadhyay, and Tiwari~\cite{usha2020evolution}, Mukhopadhyay, Ruyer-Quil, and Usha~\cite{mukhopadhyay2022modelling} have provided new insights into how two-equation models behave, particularly when attempting to apply them beyond their expected range of validity.

  Thin liquid films have numerous industrial applications, most notably in coating flows~\cite{weinstein2004coating}, heat and mass transfer~\cite{wayner1991effect}, thin-film thermoelectric cooling~\cite{darabi2001electrohydrodynamic}, as well as ice accretion prevention on aircraft surfaces~\cite{moore2017ice}. These applications require controlling the interface to a specific shape, whether that be flat for coating flows or highly corrugated for heat transfer. There are almost as many physical input mechanisms as there are applications, and there is an extensive body of literature dealing with the effects that these have on stability and the critical Reynolds number above which unstable modes exist. These include heating and cooling of the fluid~\cite{bankoff1971stability}, electric~\cite{tseluiko2006wave} or magnetic~\cite{amaouche2013hydromagnetic} fields, porous~\cite{ogden2011gravity} or deformable~\cite{gaurav2007stability} walls, and many more. Here, as shown in \cref{fig:geometry}, we focus on blowing and suction through the base via discrete actuators~\cite{thompson2016falling}, since they act over faster (and therefore more computationally accessible) timescales, and their effect on the overall flow control is greater. While continuous controls, i.e., controls applied through the entire domain (flow base), are more mathematically tractable, discrete controls applied through small holes or slots are a necessary step towards real applications, since applying blowing and suction controls throughout the whole domain is unfeasible.

  In the last two decades, the field has progressed from studying the effects that static, predetermined alterations to the system (such as fixed heating patterns or corrugated baseplates) can produce towards feedback control, where information from the interface is used to update control inputs as the film evolves in time. Armaou and Christofides~\cite{armaou2000feedback,christofides2000global} and more recently Gomes, Papageorgiou, and Pavliotis~\cite{gomes2017stabilizing} and Gomes et al.~\cite{gomes2017controlling,gomes2015controlling}, studied the simplest thin-film model, the Kuramoto-Sivashinsky (KS) equation, from both analytic and computational perspectives, and successfully applied methods from linear control theory based on~\cite{zabczyk2020mathematical}. Thompson et al~\cite{thompson2016stabilising} further extended these results to the Benney~\cite{benney1966long} and weighted-residual~\cite{ruyer2000improved} equations using a family of linear-quadratic regulator (LQR) methods. However, these long-wave models remain one stage removed from the physical system that they approximate. Unfortunately, the full two-phase Navier-Stokes system is too complex and nonlinear to directly extend the previous work on long-wave models. Cimpeanu, Gomes, and Papageorgiou~\cite{cimpeanu2021active} first analysed aspects of model reduction applicability, and delineated the discrepancies between thin-film modelling and direct numerical simulation (DNS) approaches. They considered a scenario in which actuator input is simply proportional to the observed interfacial deviation at its position, which previous numerical evidence suggested could successfully stabilise the unperturbed (flat) state of the Benney and weighted-residual systems~\cite{thompson2016stabilising}, the so-called Nusselt solution. They then showed how this model information can be interpreted and transferred into a more accurate simulation framework to successfully stabilise the full Navier-Stokes system via DNS. Nevertheless, a rigorous optimal control approach capable of surpassing the limitations of proportional control setups remains a challenge yet to be addressed. An alternative approach which incorporated model predictive control (MPC) was proposed by Wray, Cimpeanu and Gomes \cite{wray2022electrostatic} in the context of electrostatic actuation. This includes designing optimal controls for reduced-order models and enhanced interaction between model and simulation techniques, where any re-initialisation of the optimal control problem uses readings from a direct numerical simulation of the full problem. The computational cost of this framework may however prove prohibitive in real-world contexts.

  In this work we aim to close the gap between robust controls designed for long-wave models and more physically relevant systems such as those governed by the two-phase Navier-Stokes equations, ultimately bringing real-world applications closer within reach. The paper is structured as follows. In \cref{sec:governing_equations} we begin with a description of the models that make up the two-tier hierarchical structure of the control framework: the two-phase Navier-Stokes equations as the target system and either the Benney or weighted-residual model as the control system. In \cref{sec:control_methodology} we outline the control method: a set of discrete actuators injecting and removing fluid at the base of the film. We further simplify these reduced-order models to a system of linear partial differential equations (PDEs), and finally to a finite set of ODEs, which permits the use of a linear-quadratic regulator, an established control-theoretical technique. We then demonstrate for the first time that, by applying this control strategy to the full Benney, weighted-residual, and Navier-Stokes systems, there is good agreement between the linear predictions and a series of nonlinear numerical experiments. Finally, in \cref{sec:analysis} we illustrate that this agreement spans a large region of the parameter space corresponding to physically relevant fluids. Furthermore, our stability analysis results indicate that the control performance significantly exceeds the range of validity of the underpinning model assumptions in certain regions of the explored parameter space.

\section{Governing Equations}
  \label{sec:governing_equations}
\begin{figure}[ht]
    \centering%
    \includegraphics[width=\textwidth,page=1,trim={0 90 0 80},clip]{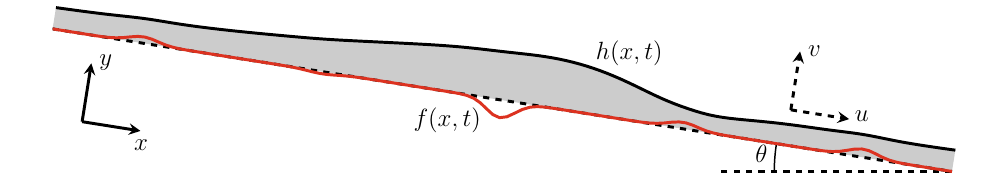}%
    \caption{Diagrammatic representation of a falling liquid film under gravity with basal forcing \(f\). Controls are applied through the wall at \(y=0\), and the behaviour of the interface \(h\) is governed by the fluid parameters and the inclination angle \(\theta\).}%
    \label{fig:geometry}%
  \end{figure}

  We consider a thin film of fluid flowing down a plane tilted at an angle \(\theta\) from the horizontal, as shown in \cref{fig:geometry}. We restrict ourselves to 2D flows, largely to make the problem more computationally tractable. Nevertheless, this setup exhibits a highly nontrivial and physically rich behaviour, while also capturing the initial wave development stages before cross-flow effects begin to appear in 3D contexts~\cite{kalliadasis2011falling}. In many control scenarios manipulating the dynamics of these early stages is the key objective (and only realisable strategy) of the control framework before highly nonlinear and often undesirable flow features arise.
  
  We use a coordinate system rotated with the plane, where \(x\) points downstream and \(y\) is the perpendicular distance to the wall. There is a free surface at the upper interface of the fluid at \(y=h(x,t)\) where the fluid and gas meet. We inject and remove fluid through the rigid lower wall regions as dictated by our resulting control strategy, with no-slip and impermeability conditions governing the remaining uncontrolled boundary. Finally, we consider periodic boundary conditions in the \(x\)-direction; while an experiment would be realised on an open domain with inflow and outflow, the speed with which a wave fully develops after the inlet~\cite{denner2018solitary} means that for sufficiently large domains -- which we consider here -- periodic boundary conditions provide a reasonable approximation. Furthermore, periodicity allows us to perform the analysis performed in section \ref{sec:analysis} below, which would be more challenging, if not impossible, to undertake given the inability to compute eigenvalues of the problem explicitly with open boundaries.

  The problem is governed by the conservation of mass and momentum in both the liquid film and the gas above it, coupled at the interface. Typically, the large density and viscosity ratios between the two media mean that we can consider the gas region to be hydrodynamically passive, and can ignore the flow in the gas, and model the liquid film alone. The fluid flow is governed by the acceleration due to gravity \(g\), the inclination angle \(\theta\), and the physical properties of the liquid phase: constant density \(\rho\), viscosity \(\mu\) and surface tension coefficient \(\gamma\). A full list of physically-relevant parameters and the values used in this investigation can be found in \cref{sec:parameter_values}.

  For a liquid film with mean height \(h_s\), the uncontrolled system admits a uniform solution known as the Nusselt solution~\cite{nusselt1923warm}, where \(h(x,t)=h_s\), which has a parabolic velocity profile with surface velocity \(U_s = \frac{\rho g h_s^2\sin\theta}{2\mu}\). We then nondimensionalise the problem based on the length scale \(h_s\), velocity scale \(U_s\) and pressure scale \(\frac{\mu U_s}{h_s}\), defining the Reynolds and capillary numbers
  \begin{equation}
    \Rey = \frac{\rho U_s h_s}{\mu}, \qquad \Ca = \frac{\mu U_s}{\gamma},
  \end{equation}
  which measure the relative importance of inertial and viscous terms, and of gravity and surface tension, respectively.

\subsection{Navier-Stokes equations}
  \label{sub:navier_stokes_equations}
  The full liquid film flow is governed by the 2D Navier-Stokes equations, which are solved for velocity \(\vec{u}(x,y,t) = (u,v)\) and pressure \(p(x,y,t)\) under the action of external forces. The governing (nondimensionalised) internal momentum equations are
  \begin{align}
    \label{eqn:momentum1}\Rey (u_t + uu_x + vu_y) &= -p_x + 2 + u_{xx} + u_{yy}, \\
    \label{eqn:momentum2}\Rey (v_t + uv_x + vv_y) &= -p_y - 2\cot\theta + v_{xx} + v_{yy},
  \end{align}
  and the continuity equation reads
  \begin{equation}
    \label{eqn:continuity}u_x + v_y = 0.
  \end{equation}

  The system is completed by its boundary conditions: periodic boundaries in the \(x\)-direction, no-slip and fluid injection/removal at the wall,
  \begin{equation}
    \label{eqn:wall}u = 0, \qquad v = f(x,t),
  \end{equation}
  the nonlinear dynamic stress balance (or momentum jump) at the interface, \(y=h(x,t)\),
  \begin{align}
    \label{eqn:stress1}(v_x + u_y)(1-h_x^2) + 2h_x (v_y  -u_x) &= 0, \\
    \label{eqn:stress2}p - \frac{2}{1+h_x^2}(v_y + u_x h_x^2 - h_x(v_x + u_y)) &= -\frac{1}{\Ca}\frac{h_{xx}}{(1+h_x^2)^{3/2}},
  \end{align}
  and finally the kinematic boundary condition
  \begin{equation}
    \label{eqn:kinematic}h_t = v-uh_x.
  \end{equation}

  In lieu of physical experiments, we perform computational analogues by simulating the Navier-Stokes equations using a volume-of-fluid approach developed by Popinet~\cite{popinet2003}. The methodology is well known, following more than two decades of successful development and usage in the community~\cite{popinet2003,popinet2009,popinet2015}, and therefore we restrict our attention to details relevant to our particular setting in \cref{sec:numerical_simulations}.

  Defining the down-slope flux \(q(x,t)\) by integrating over the height of the film
  \begin{equation}
    q(x,t) = \int_{0}^{h} u(x,y,t) \diff y,
  \end{equation}
  we combine \cref{eqn:continuity,eqn:wall,eqn:kinematic} to obtain the 1D mass conservation equation
  \begin{equation}
    \label{eqn:conservation}h_t + q_x = f.
  \end{equation}

  We can continue to use the Navier-Stokes equations to compute \(q\), or we can use one of a number of simplified models for the flux. Here we consider the Benney~\cite{benney1966long} and weighted-residual~\cite{ruyer2000improved} equations, which are valid in the long-wave limit. By using a pair of reduced-order models we are better able to gauge the relative capabilities of both, weighing model and computational complexity against control performance. The following pair of reduced-order models are based on first-order asymptotic expansions in the long-wave parameter \(\epsilon = 1/L\) (where $L$ is the aspect ratio of the domain). In addition to the requirement that \(\epsilon \ll 1\), we make the assumption that \(\Rey = O(1)\) and \(\Ca = O(\epsilon^2)\) to retain inertial and surface tension effects, and that \(f = O(\epsilon)\) so that the magnitude of the imposed control is comparable to the perturbed flow.

\subsection{Benney equation}
  \label{sub:benney_equation}
  The first choice of model for the downstream flux is the Benney system~\cite{benney1966long}, which was extended to include the effects of \(O(\epsilon)\) controls by Thompson, Tseluiko, and Papageorgiou~\cite{thompson2016falling}. This enslaves the flux to the interfacial height \(h\) via
  \begin{equation}
    \label{eqn:benney}q(x,t) = \frac{h^3}{3}\left( 2 - 2 h_x \cot\theta + \frac{h_{xxx}}{\Ca} \right) + \Rey \left( \frac{8h^6h_x}{15} - \frac{2h^4f}{3} \right),
  \end{equation}
  resulting in a single equation for the evolution of the interface when coupled to \cref{eqn:conservation}. The system above is a significant improvement over equations such as the KS equation -- the simplest nonlinear model of thin film flows~\cite{sivashinsky1980irregular}. While capturing some important aspects of the chaotic behaviour of falling liquid films such as travelling waves, the KS equation is more useful as a paradigmatic example in a dynamical systems sense than as a predictive model outside of a very restricted region of the parameter space. However, the Benney system exhibits undesirable behaviours such as unphysical finite-time blow-up outside of a narrow range of parameters corresponding to low Reynolds numbers~\cite{pumir1983solitary}, as is demonstrated in \cref{fig:interface}.

\subsection{Weighted-residual system}
  \label{sub:weighted_residual_system}
\begin{figure}[ht]
    \centering%
    \includegraphics[width=\textwidth,page=4,trim={30 15 20 15},clip]{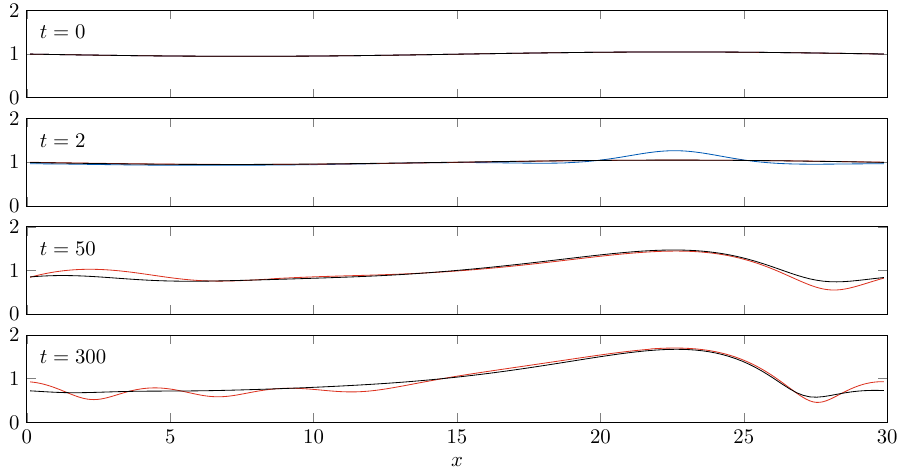}%
    \caption{Evolution of interfacial heights \(h\) for Navier-Stokes (black), weighted-residual (red), and Benney (blue) systems, with peaks shifted to \(3L/4\). Here, the parameters used are \(\Rey=10\), \(\Ca=0.05\), \(\theta=\pi/3\). The Benney equation blows up shortly after \(t=2\), but the weighted-residual and Navier-Stokes interfaces have very similar structures aside from some spurious oscillations in the weighted-residual case, which are observable at $t=300$ above.}%
    \label{fig:interface}%
  \end{figure}

 To overcome the unrealistic behaviour described above, Ruyer-Quil and Manneville~\cite{ruyer2000improved,ruyer2002further} proposed an improved weighted-residual methodology based on approximating \(u\) by a truncated sum of basis functions satisfying no-slip boundary conditions at the wall~\eqref{eqn:wall} and zero tangential stress at the interface~\eqref{eqn:stress1}. Here we use the first-order truncation, which matches well with the second-order truncation up to \(\Rey \approx 5\)~\cite{ruyer2000improved}, a significant improvement over previous models~\cite{salamon1994traveling}. When combined with the basal forcing this results in a separate evolution equation for the flux~\cite{thompson2016falling}
  \begin{equation}
    \label{eqn:wr}\frac{2\Rey}{5}h^2q_t + q = \frac{h^3}{3}\left( 2 - 2 h_x \cot\theta + \frac{h_{xxx}}{\Ca} \right) + \Rey \left( \frac{18q^2h_x}{35} - \frac{34hqq_x}{35} + \frac{hqf}{5} \right),
  \end{equation}
  which together with~\cref{eqn:conservation} forms a system of two PDEs for the height \(h(x,t)\) and the flux \(q(x,t)\). \Cref{eqn:wr} better captures many features of the full Navier-Stokes film, such as spontaneous back-flow~\cite{oron2009numerical}, and \cref{fig:interface} illustrates how it provides a good match for the interfacial shape even at moderate Reynolds and capillary numbers. However, it does overestimate the amplitude of the capillary ripples, as observed by Ruyer-Quil and Manneville~\cite{ruyer2000improved} for the weighted-residuals system and also for other first-order models~\cite{alekseenko1994wave}. Unlike the Benney equation, it also does not exhibit unphysical finite-time blow-up, although it too diverges from the Navier-Stokes model at moderate Reynolds numbers, \(\Rey \approx 10\)~\cite{ruyer2000improved}.
  
\section{Control methodology}
  \label{sec:control_methodology}
  We focus on controlling the interface towards the Nusselt solution, which under our nondimensionalisation is the uniform film \(h(x,t)=1\). All the controls we consider are a class of time-dependent controls known as feedback controls, which we introduce here. Take a controlled quantity \(x\) governed by the system
  \begin{equation}
    x_t = \mathcal{A}x + \mathcal{B}u, \qquad y = \mathcal{C}x,
  \end{equation}
  where \(u\) is the control, \(y\) is some observation of the system, and \(\mathcal{A}\), \(\mathcal{B}\), and \(\mathcal{C}\) are arbitrary operators describing the uncontrolled system dynamics, the control actuation mechanism, and the observations respectively. In the case of feedback controls, we have the restriction that \(u = \mathcal{K}y\) for some operator \(\mathcal{K}\), so that the system can be written in closed-loop form
  \begin{equation}
    \label{eqn:closed_loop}x_t = (\mathcal{A} + \mathcal{B}\mathcal{K}\mathcal{C})x.
  \end{equation}
  An overview of some important control theory definitions is provided in \cref{sec:control_theory_fundamentals}.

\begin{figure}[ht]
    \centering%
    \includegraphics[width=\textwidth,page=2,trim={5 20 5 20},clip]{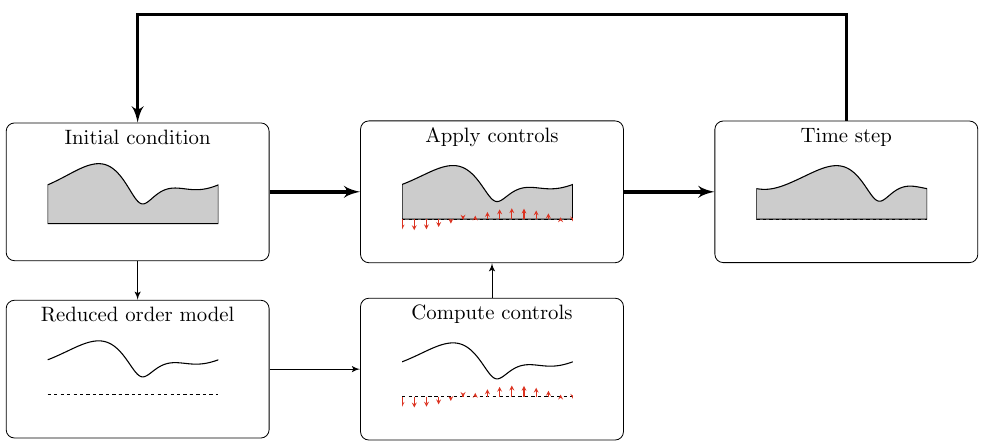}%
    \caption{Multi-layer control methodology for the control of Navier-Stokes thin liquid films. From the initial condition, we treat the interface as though it were described by the chosen reduced-order model and generate the feedback control accordingly. We then apply this control to the full model and time step forward to repeat the process.}%
    \label{fig:framework}%
  \end{figure}

  In the case of falling liquid films, the full system \cref{eqn:momentum1,eqn:momentum2,eqn:continuity,eqn:wall,eqn:stress1,eqn:stress2,eqn:kinematic} is too complex for standard (linear) control-theoretical techniques to be tractable. Instead, we design feedback controls for the reduced-order models presented in \cref{sub:benney_equation,sub:weighted_residual_system} and apply them to the full Navier-Stokes system by passing the Navier-Stokes interfacial height to the feedback control scheme. The full framework is pictured in \cref{fig:framework}.

  Given the difficulties in observing both the height~\cite{liu1993onset,vlachogiannis2002experiments} and flux~\cite{heining2013flow} of falling liquid films, it is understandable that in our case we might wish to express our control \(f\) as a function of time only, known as \rnew{open-loop or} offline control (as can be done when generating controls for the KS equation for instance -- see the work of Gao \cite{gao2016optimal}). Unfortunately, \rnew{this is known to be impractical. A}s shown by Cimpeanu, Gomes, and Papageorgiou~\cite{cimpeanu2021active}, although such hierarchical controls show promising initial dampening, they invariably fail over longer timescales, as the Navier-Stokes model eventually diverges from the chosen long-wave model, and the action of the control no longer affects the intended state, with the Navier-Stokes system ultimately converging back to its uncontrolled behaviour.

\subsection{Control actuation mechanism}
  \label{sub:control_mechanism}
  Although the control term \(f\) in \cref{eqn:conservation,eqn:benney,eqn:wr} is general, in this study we restrict ourselves to controls taking the form of Dirac-delta distributions injecting or removing fluid at the wall at a finite number of locations \(x_1, \ldots, x_M\) -- which is more experimentally achievable than continuous controls, i.e. controls applied everywhere in the domain. Furthermore, both due to realistic considerations and computational restrictions imposed by the direct numerical simulation setup, we must approximate these point sources by finite regions which we select to be smooth, periodic functions (shown in \cref{fig:gain_matrix}), of the type
  \begin{equation}
    d(x) = A \exp\left[ \frac{\cos(2\pi x/L)-1}{\omega^2} \right],
  \end{equation}
  where \(\omega\) controls the width of the function, and \(A\) is chosen so that \(\int_{0}^{L} d(x) \diff x = 1\). More refined discretisations support smaller values of \(\omega\), and \(d(x)\to\delta(x)\) as \(\omega \to 0\). The basal forcing term \(f\) is thus
  \begin{equation}
    f(x, t) = \sum_{i=1}^{M} u_i(t) d(x-x_i),
  \end{equation}
  where \(u_i(t)\) are the individual, time-dependent, control amplitudes. Despite the practical difficulties of obtaining full observations, the main goal of this investigation is to test the feasibility of the \emph{control} methodology, and so for the moment we assume we are able to observe the full interface \(h(x)\). This means that \(\mathcal{C}\), the observation operator from \cref{eqn:closed_loop}, is the identity. We will address the \emph{observability} of the problem, as well as issues introduced by noisy or partial observations of the interface, in future work.
  
  Note that, since the domain has periodic boundaries on the left and right, the problem is translationally invariant, and so \(x_1, \ldots, x_M\) should be evenly placed along the base. For 3D and non-periodic flows, the optimal placement of the actuators is a nontrivial problem~\cite{tomlin2019point}.

  Finally, we introduce a cost functional to compare different control strategies, taking into account the 2-norms of the deviation from the target state and penalising the use of the controls. We thus define the total cost of a control by
  \begin{equation}
    \label{eqn:cost}\kappa = \int_{0}^{\infty} \int_{0}^{L}  \beta \hat{h}(x)^2 + (1-\beta) f^2 \diff x \diff t,
  \end{equation}
  where \(\hat{h}(x) = h(x)-1\) is the deviation of the interface from the target uniform state, and the parameter \(\beta\) controls the relative importance of the interfacial deviation and the magnitude of the controls.

\subsection{Linear-Quadratic Regulator (LQR)}
  \label{sub:lqr}
  Despite the long-wave simplification to one of the two reduced-order models, control methodologies for nonlinear PDEs of the type considered here are still a rapidly developing area of active research, with the most relevant efforts by Boujo and Sellier~\cite{boujo2019pancake}, Lunz~\cite{lunz2021minimizing}, or the current authors~\cite{cimpeanu2021active,gomes2017stabilizing,thompson2016stabilising}. In spite of these recent advances, we cannot directly choose the optimal control operator \(\mathcal{K}\) with analytical methods. To make this problem tractable, we assume that the perturbation away from the Nusselt solution (\(\hat{h}=h-1=0\), \(\hat{q}=q-2/3=0\), \(\hat{f}=0\)) is small, so that we can linearise \cref{eqn:benney} to obtain the equation
  \begin{equation}
    \label{eqn:lin_benney}\hat{h}_t = \left[ -2\partial_x + \left( \frac{2\cot{\theta}}{3} - \frac{8\Rey}{15} \right)\partial_{xx} - \frac{1}{3\Ca}\partial_{xxxx} \right]\hat{h} + \left[ 1 + \frac{2\Rey}{3}\partial_x \right] \hat{f}.
  \end{equation}
  We then discretise to form a system of \(N\) ODEs,
  \begin{equation}
    \label{eqn:lqr}\frac{\diff \hat{h}}{\diff t} = J\hat{h} + \Psi u, \qquad u = K\Phi\hat{h}.
  \end{equation}
  \rnew{Importantly,  Gibson~\cite{gibson1983linear} showed that, despite relying on a linearisation and discretisation of the original system, under certain conditions the discretised feedback operator \(K\) obtained using the approach described below converges to its infinite-dimensional counterpart \(\mathcal{K}\). Although the linearised and discretised system is now far removed from the original Navier-Stokes equations, this approach is enabled by the ability of this operator to capture the undesirable instabilities sufficiently well so as to provide a basis for effective control.}
  
  Here, \(J\in\mathbb{R}^{N\times N}\) captures the system dynamics, \(\Psi\in\mathbb{R}^{N\times M}\) is the linearised actuator matrix, and \(\Phi\in\mathbb{R}^{N\times N}\) is the linearised observation matrix (which we take to be the identity). \(K\in\mathbb{R}^{M\times N}\) is the gain matrix, which is chosen to minimise the discrete cost
  \begin{equation}
    c = \int_{0}^{\infty} \hat{h}\tran U \hat{h} + u\tran V u \diff t,
  \end{equation}
  where \(U = \frac{\beta L}{N} I\in\mathbb{R}^{N\times N}\) and \(V = (1-\beta) I\in\mathbb{R}^{M\times M}\) are matrices whose entries are chosen so as to form the discrete analogue of the continuous cost \cref{eqn:cost}. The process is similar for the weighted-residual system, but with twice the system size at each stage, since there are two unknowns, $\hat{h}$ and $\hat{q}$. The linearisation of \cref{eqn:wr} results in
  \begin{align}
    \label{eqn:lin_wr1}\hat{h}_t &= -\hat{q}_x + \hat{f}, \\
    \label{eqn:lin_wr2}\hat{q}_t &= \left[ \frac{5}{\Rey} + \left( \frac{4}{7} - \frac{5\cot{\theta}}{3\Rey} \right)\partial_x + \frac{5}{6\Rey\Ca}\partial_{xxx} \right]\hat{h} - \left[\frac{5}{2\Rey} + \frac{34}{21}\partial_x \right]\hat{q} + \left[ \frac{1}{3} \right]\hat{f},
  \end{align}
  and the resulting discretised system has \(2N\) equations rather than \(N\). Finally, although we are assuming full observations of the interfacial height, Thompson et al.~\cite{thompson2016stabilising} showed that it is sufficient to use the leading order approximation \(\hat{q} = 2\hat{h}\) to remove the need to directly observe the flux, incurring a small penalty in the size of the largest eigenvalue but not fundamentally affecting stability. This is particularly important because the flux is challenging to measure in an application setting~\cite{heining2013flow}.

  This setup forms a classic problem in control theory: the linear-quadratic regulator (LQR) problem, which is a subset of a broader class of static output feedback (SOF) problems in which one can also have restricted observations (i.e. \(\text{rank}(\Phi) < N\)). Here, we provide an overview of how this class of problems is solved. For more details see~\cite{johnson1970design,syrmos1997static}.

  For the discretised linear control system \cref{eqn:lqr}, we write the cost as
  \begin{equation}
    \label{eqn:cost_sub}c = \int_{0}^{\infty} \hat{h}\tran U \hat{h} + u\tran V u \diff t \\
    = \int_{0}^{\infty} \hat{h}\tran (U + \Phi\tran K\tran V K \Phi) \hat{h} \diff t,
  \end{equation}
  where \(U, V\) are assumed to be symmetric positive definite matrices.

  If we suppose there exists a symmetric, positive semi-definite matrix \(P\) such that
  \begin{equation}
    \label{eqn:defP}\frac{\diff}{\diff t} (\hat{h}\tran P \hat{h}) = -\hat{h}\tran (U + \Phi \tran K\tran V K \Phi) \hat{h},
  \end{equation}
  then, as long as the controlled system matrix \(A = J+\Psi K \Phi\) is asymptotically stable, i.e., all its eigenvalues have negative real part, we can write \cref{eqn:cost_sub} as
  \begin{equation}
  \begin{split}
    c &= \hat{h}(0)\tran P \hat{h}(0) - \lim_{t\to\infty} \hat{h}(t)\tran P \hat{h}(t) \\
    &= \hat{h}(0)\tran P \hat{h}(0).
  \end{split}
  \end{equation}
  By expanding out the left hand side of \cref{eqn:defP} and observing that this is true for all initial conditions \(\hat{h}(0)\in\mathbb{R}^N\), we have
  \begin{equation}
    \label{eqn:constraint}A\tran P + P A + U + \Phi\tran K\tran V K \Phi = 0.
  \end{equation}
  This further implies that the choice of \(P\) is independent of the initial condition \(\hat{h}(0)\), and so
  \begin{equation}
    \label{eqn:cost_trace}c = \tr(PX),
  \end{equation}
  where \(X = \hat{h}(0)\hat{h}(0)\tran\). Since we wish to choose an optimal \(K\) for all initial conditions, we set \(X = \mathbb{E}[\hat{h}(0)\hat{h}(0)\tran] = I\), the identity matrix, as we assume all initial perturbations \(\hat{h}(0)\) are equally likely.

  The problem thus becomes equivalent to selecting \(K\) to minimise \cref{eqn:cost_trace} subject to the constraint \cref{eqn:constraint}. This can be solved via Lagrange multipliers. Defining the symmetric matrix of Lagrange multipliers \(S\), we then have the resulting Hamiltonian
  \begin{equation}
    H = \tr(PI) + \tr((A\tran P + P A + U + \Phi\tran K\tran V K \Phi)S).
  \end{equation}
  By setting \(\partial_S H = \partial_P H = \partial_K H = 0\) we have the conditions for the solution to the SOF problem:
  \begin{align}
    \label{eqn:cond1}0 &= A\tran P + P A + U + \Phi\tran K\tran V K \Phi , \\
    \label{eqn:cond2}0 &= A S + S A\tran + I, \\
    \label{eqn:cond3a}0 &= VK\Phi S\Phi\tran + \Psi\tran PS\Phi\tran.
  \end{align}
  The final condition can be more usefully written as
  \begin{equation}
    \label{eqn:cond3}K = -V^{-1}\Psi\tran P S\Phi\tran (\Phi S\Phi\tran)^{-1}.
  \end{equation}

  \Cref{eqn:cond1,eqn:cond2,eqn:cond3a} cannot be solved directly, and so an iterative procedure must be used. However, in the special case of the LQR problem where we have \(\Phi = I\), \rchange{ we may }{the expression \cref{eqn:cond3} for \(K\) no longer depends on \(S\). We may thus} discard \cref{eqn:cond2} and rewrite \cref{eqn:cond1,eqn:cond3} as
  \begin{align}
    \label{eqn:CARE}0 &= J\tran P + P J + U -P\Psi V^{-1}\Psi\tran P, \\
    \label{eqn:gain_matrix}K &= -V^{-1}\Psi\tran P.
  \end{align}

  \Cref{eqn:CARE}, which is known as the continuous algebraic Riccati equation (CARE), can be solved for \(P\) directly, and then used to compute \(K\). The structure of the matrices \(J\), \(U\), \(V\), and \(\Psi\) -- with \(U\) and \(V\) diagonal, \(J\) periodic banded and \(\Psi\) having translational symmetries --  means that the specific CARE for this problem is typically well-conditioned. Thus we can make use of the classical eigenvector approach described by MacFarlane~\cite{macfarlane1963eigenvector}, Potter~\cite{potter1966matrix} and Vaughan~\cite{vaughan1970nonrecursive}. Alternatively, Schur~\cite{laub1979schur} and generalised eigenvector~\cite{arnold1984generalized} approaches may offer improved numerical stability for larger systems (which would be encountered in 3D) and more unstable regimes (where some of the interim matrices used in the classical method become singular or near-singular).

  We note that equations~\cref{eqn:cost_trace,eqn:gain_matrix}\rnew{, alongside the lack of a cross-term in the cost expression \cref{eqn:cost_sub}}, illustrate why the single parameter \(\beta\) is sufficient to fully explore the cost-space with regards to \(K\): if we \rchange{instead introduce a pair of control parameters $\alpha$ and $\beta$}{introduce an additional control parameter \(\alpha\)} so that
  \begin{equation}
    U^\prime = \alpha U = \frac{\alpha\beta L}{N} I, \qquad V^\prime = \alpha V = \alpha(1-\beta) I,
  \end{equation}
  we can set the entries of \(U^\prime\) and \(V^\prime\) independently. The cost is then
  \begin{equation}
    c^\prime = \alpha c = \frac{1}{2} \tr(\alpha P X) = \frac{1}{2} \tr(P^\prime X).
  \end{equation}
  Carrying the new cost matrices through to \cref{eqn:gain_matrix} we have
  \begin{equation}
    K^\prime = -(V^\prime)^{-1}\Psi\tran P^\prime = -(\alpha V)^{-1}\Psi\tran \alpha P = K.
  \end{equation}
  The above result indicates that scaling the cost makes no difference to the optimal \(K\), and so a single parameter describing the ratio of significance of the two components is sufficient.\rold{ Gibson~\cite{gibson1983linear} showed that, under certain conditions, the discretised feedback operator $K$ does converge to its infinite-dimensional counterpart $\mathcal{K}$.}

\begin{figure}[ht]
    \centering%
    \includegraphics[width=\textwidth,page=3,trim={10 20 10 20},clip]{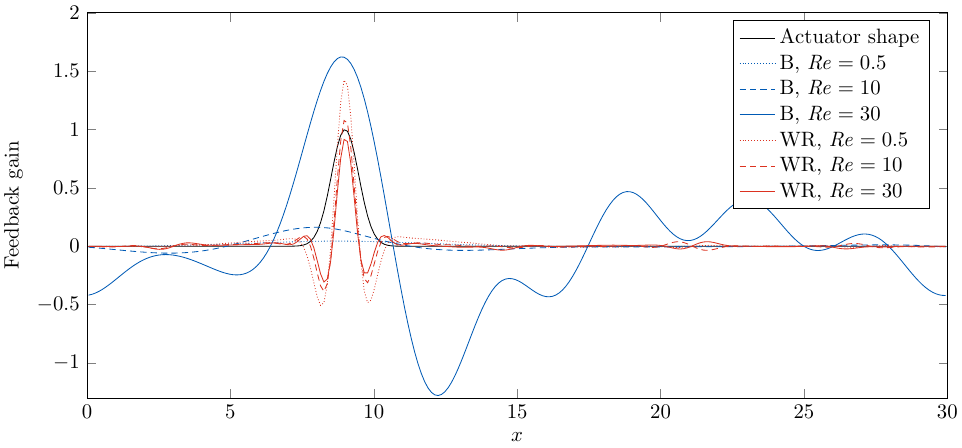}%
    \caption{The second row of the gain matrix computed using either the Benney equation (in blue) or weighted-residual system (in red) as the reduced-order model, as \(\Rey\) varies and \(Ca\) is fixed at 0.05. The gains are shown alongside the corresponding actuator (in black). Although the weighted-residual gains remain clustered around the actuator, the Benney gains have significant nonlocal contributions. \rnew{The strong correlation between actuator location and feedback gain suggests that the discretisation provides a good approximation to the continuous analogue.}}%
    \label{fig:gain_matrix}%
  \end{figure}

  Once the optimal gain matrix \(K\) has been computed, we can calculate the \(m^{\textrm{th}}\) actuator amplitude as \(u_m = K_m \cdot \hat{h}(t)\), where \(K_m\) is the \(m^{\textrm{th}}\) row of \(K\) and \(\cdot\) denotes the inner product. This means that \(K_{m,i}\) can be interpreted as describing the importance of the \(i^{\textrm{th}}\) entry of \(\hat{h}\) to \(u_m\). As can be seen in \cref{fig:gain_matrix}, this allows us to examine the rows of \(K\) to develop an understanding of how the controls operate. The weighted-residual gains are tightly clustered around the location of the actuator across a wide range of Reynolds numbers, with minimal up- and down-stream contributions. By contrast, the Benney gains are much broader and depend more strongly on the interfacial shape away from the actuator location. They also are much more sensitive to the Reynolds number (it is worth noting that, for \(\Rey = 30\), the Benney-derived controls fail to stabilise the Navier-Stokes system).

 With a method to compute the gain matrix for the two reduced-order models we are now well-positioned to deploy the methodology described in \cref{fig:framework} and direct it towards the modelled physical system of interest.

\subsection{Preliminary results}
  \label{sub:results}
\begin{figure}[ht]
    \centering%
    \includegraphics[width=\textwidth,page=8,trim={10 15 5 15},clip]{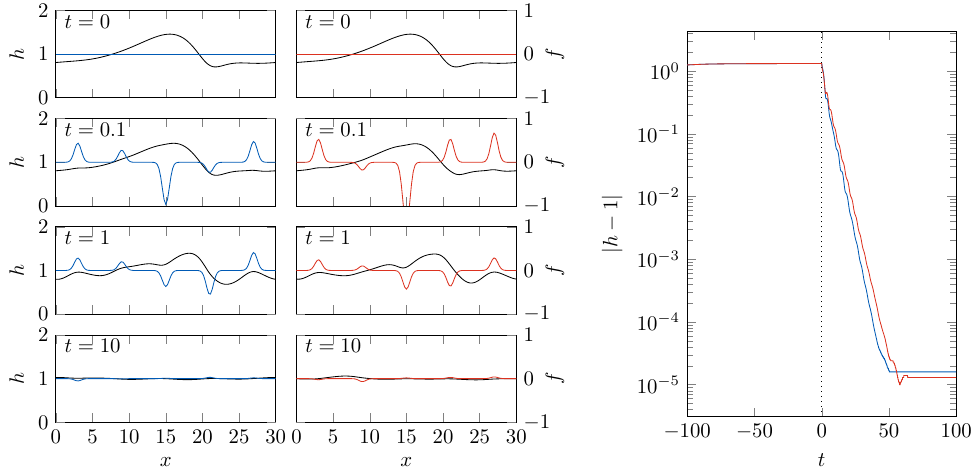}%
    \caption{Interfacial shapes before and after the controls are switched on: Benney equation derived controls in blue (left), weighted-residual derived controls in red (centre). A travelling wave is allowed to develop until \(t=0\), when the controls are activated. Both controls successfully damp out the perturbation, with the control amplitudes decreasing in proportion to \(\abs{h-1}\). We note that, although similar, the controls are not identical -- see the second and fourth rows in particular. In both cases (Benney in blue, weighted-residual in red) the 2-norm of the deviation of the interface from the target state decays exponentially (right). After \(t\approx 50\) the deviation is small enough that machine precision interferes with computing the deviation. In these simulations we used \(\Rey = 5\), \(\Ca = 0.05\).}%
    \label{fig:controlled_interfaces}%
  \end{figure}

  Previous work by Thompson et al.~\cite{thompson2016stabilising} confirmed that LQR controls with full observations are able to stabilise both the Benney and weighted-residual systems\rnew{ (as given by equations \cref{eqn:conservation} with \cref{eqn:benney} or \cref{eqn:wr}, respectively)}. The same authors also found that the Benney controls stabilise the weighted-residual model. In \cref{fig:controlled_interfaces} we can see that for a similar parameter regime (\(\Rey = 5\), \(\Ca = 0.05\) -- selected such that \(\Rey\) is not so high so as to make numerical simulation difficult, and \(\Ca\) is large enough that surface tension alone cannot stabilise the liquid film with properties experimentally aligned with a relatively thick and viscous oil flow), these controls can be extended to the Navier-Stokes system, where we achieve similar results. \Cref{fig:controlled_interfaces} shows how the interface is allowed to develop from a small sinusoidal perturbation into a travelling wave, before the application of controls at \(t=0\). Representative interfacial snapshots are pictured in \cref{fig:controlled_interfaces}. The interfacial deviation then decays exponentially, suggesting that the use of linear models to design the gain matrix is appropriate in both cases.

\section{Stability analysis}
  \label{sec:analysis}
  It is encouraging to see that we can control the film in the specific setting of \cref{fig:controlled_interfaces}, but a better aim is to predict the stabilisability of the system given the flow parameters \(\Rey\), \(\Ca\), \(\theta\) and number of controls \(M\). Since we lack a closed-form expression for either the continuous control \(f(h)\), or its discrete counterpart \(\Psi K \hat{h}\), we cannot directly estimate the stability properties of the controlled system. However, we can predict the damping rate by finding the largest eigenvalue \(\lambda^*\) of the controlled system matrix \(A = J+\Psi K \Phi\) and compare that to rates fitted to the data produced in our numerical simulations. \rnew{For a given point in the parameter space, we initialise our experiment from the steady travelling wave solution, but, as illustrated in \cref{sec:initial_conditions}, we observe the same damping behaviour regardless of the initial condition}.

  \begin{figure}[ht]
    \centering%
    \includegraphics[width=\textwidth,page=5,trim={10 20 10 20},clip]{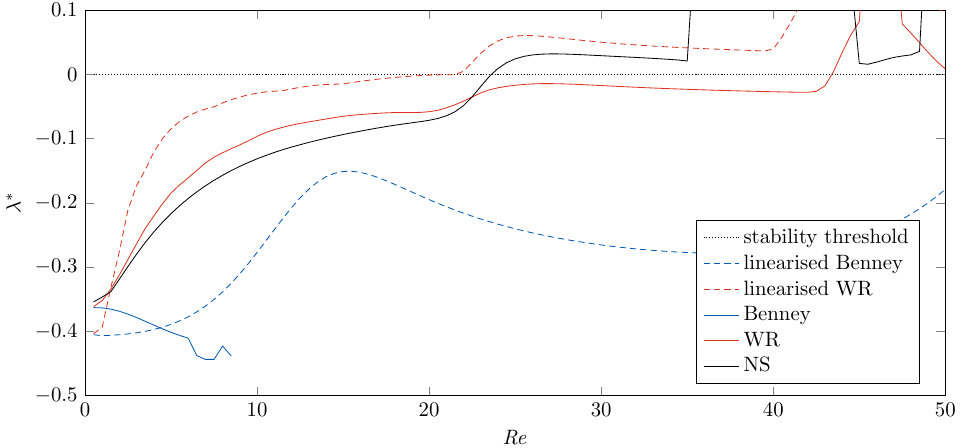}%
    \caption{Comparison of (fitted) damping rates for Benney-derived LQR control applied to Benney (blue), weighted-residual (red), and Navier-Stokes (black) systems (all solid) to the predictions from the linearised systems of ODEs. Here, we used \(M=5\) controls and \(\Ca=0.05\). For all three systems, the numerical models break down at sufficiently large \(\Rey\).}%
    \label{fig:rates_benney}%
  \end{figure}

  From \cref{fig:rates_benney} we observe that the Benney-derived controls directly stabilise the Benney and weighted-residual systems (in a similar setup to that used by Thompson et al.~\cite{thompson2016stabilising}) over a wide range of Reynolds numbers, and that their ability to stabilise towards the uniform film extends to the hierarchical controls applied to the Navier-Stokes film. We note that the weighted-residual and Navier-Stokes systems are stabilised even above the stability threshold, after which the linearised weighted-residual model predicts that five actuators are not sufficient to stabilise the uniform state.

  All three models display unphysical blow-up at sufficiently large Reynolds numbers. Although this is expected behaviour in the case of the Benney film~\cite{nakaya1975long}, in the case of the weighted-residual and Navier-Stokes models this is attributed to the eventual breakdown of the controls as the Benney model finally loses the last of its predictive capacity at larger values of \(\Rey\).

  \begin{figure}[ht]
    \centering%
    \includegraphics[width=\textwidth,page=6,trim={10 20 10 20},clip]{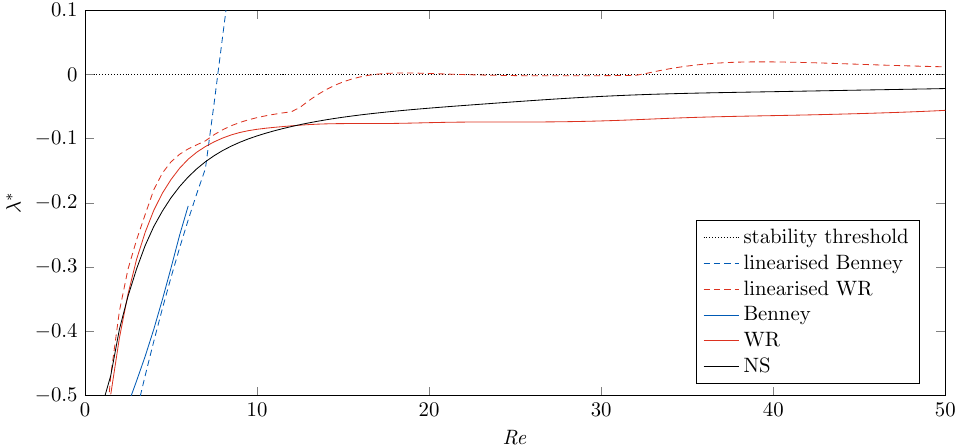}%
    \caption{Comparison of (fitted) damping rates for weighted-residual-derived LQR control applied to Benney (blue), weighted-residual (red), and Navier-Stokes (black) systems (all solid) to the predictions from the linearised systems of ODEs. Here, we used \(M=5\) controls \(\Ca=0.05\).}%
    \label{fig:rates_wr}%
  \end{figure}

  While the Benney-derived control rules stabilise all three models (at least for small-to-moderate Reynolds numbers), the weighted-residual derived controls fail to stabilise the Benney equation for \(\Rey > 7\), in agreement with the linear predictions given by the eigenvalues of \(A\). The weighted-residual and Navier-Stokes models have reasonable agreement with the linear damping rates but remain stabilisable even at \(\Rey = 50\), when the linear system is not.

  In order to make analytical progress, we turn to an equivalent way to produce the gain matrix \(K\), where we first convert \cref{eqn:lqr} to Fourier space (so \(\tilde{h} = \mathcal{F}\hat{h}\), where \(\mathcal{F}\) is the Fourier transform). We can then reorder the wavenumbers to separate stable and unstable modes:
  \begin{equation}
    \frac{\diff \tilde{h}}{\diff t} = \tilde{J}\tilde{h} + \tilde{\Psi}\tilde{K}\tilde{h} =
    \matrix{\tilde{J}_u & 0 \\ 0 & \tilde{J}_s}\tilde{h} + \matrix{\tilde{\Psi}_u \\ \tilde{\Psi}_s}\tilde{K}\tilde{h}.
  \end{equation}
  Concentrating on the unstable modes more explicitly, i.e., 
  \begin{equation}\label{eqn:FourierSplit}
  \begin{split}
    \frac{\diff}{\diff t}\matrix{\tilde{h}_u \\ \tilde{h}_s} &=  \matrix{\tilde{J}_u & 0 \\ 0 & \tilde{J}_s}\matrix{\tilde{h}_u \\ \tilde{h}_s} + \matrix{\tilde{\Psi}_u \\ \tilde{\Psi}_s}\tilde{K}\matrix{\tilde{h}_u \\ \tilde{h}_s} \\
    &= \matrix{\tilde{J}_u + \tilde{\Psi}_u \tilde{K}_u & 0 \\ \tilde{\Psi}_s \tilde{K}_s & \tilde{J}_s}\matrix{\tilde{h}_u \\ \tilde{h}_s},
  \end{split}
  \end{equation}
  we find that since the matrix on the right-hand side of \eqref{eqn:FourierSplit} is block lower triangular, the controls leave the eigenvalues of the stable modes unchanged, and so they remain stable. We thus reduce the control problem to
  \begin{equation}
    \frac{\diff \tilde{h}_u}{\diff t} = \tilde{J}_u \tilde{h}_u + \tilde{\Psi}_u \tilde{K}_u \tilde{h}_u.
  \end{equation}
  By solving the problem in Fourier space it is clear that -- for the purely linear case at least -- we should expect that \(M\) actuators would be sufficient to control any system satisfying \(M \ge \text{rank}(\tilde{J}_u)\). This would amount to  one control per unstable mode plus one more to satisfy conservation of mass, as pointed out by Armaou and Christofides~\cite{christofides2000global}.

  The rank of the unstable Jacobian \(\tilde{J}_u\) corresponds to the number of unstable modes of the linearised system (\cref{eqn:lin_benney} or \cref{eqn:lin_wr1,eqn:lin_wr2}). We compute this rank for a perturbation with wavenumber \(k\), where the linearised Benney equation has a single eigenvalue
  \begin{equation}
    \label{eqn:benney_eval}\lambda = -2ik + \left( \frac{8\Rey}{15} - \frac{2}{3}\cot\theta - \frac{1}{3\Ca}k^2 \right)k^2,
  \end{equation}
  and the weighted-residual system has a pair of eigenvalues that solve the quadratic equation
  \begin{equation}
    \label{eqn:wr_eval}\lambda^2 + \left( \frac{5}{2\Rey} + \frac{34}{21}ik \right)\lambda + \left( \frac{5}{\Rey}ik - \left[ \frac{4}{7} - \frac{5\cot\theta}{3\Rey} \right]k^2 + \frac{5}{6\Rey\Ca}k^4 \right) = 0.
  \end{equation}
  Setting the real part \(\Re(\lambda) = 0\) we can solve for the critical wavenumber \(k_0\) (the boundary between stable and unstable unimodal systems). For both \cref{eqn:benney_eval} and \cref{eqn:wr_eval}, this is
  \begin{equation}
    k_0 = \pm\sqrt{\Ca\left(\frac{8}{5}\Rey-2\cot\theta\right)}.
  \end{equation}
  After rescaling to account for \(L\neq 2\pi\), this expression admits a single zero eigenmode and pairs of positive and negative modes with \(k<k_0\), resulting in the number of unstable modes being
  \begin{equation}
    \label{eqn:unstable_modes}n_{\text{u}} = 1 + 2\left\lfloor \frac{L}{2\pi} k_0 \right\rfloor = 1 + 2\left\lfloor \frac{L}{2\pi}\sqrt{\Ca\left(\frac{8}{5}\Rey-2\cot\theta\right)} \right\rfloor.
  \end{equation}

\begin{figure}[ht]
    \centering%
    \includegraphics[width=\textwidth,page=7,trim={20 20 20 20},clip]{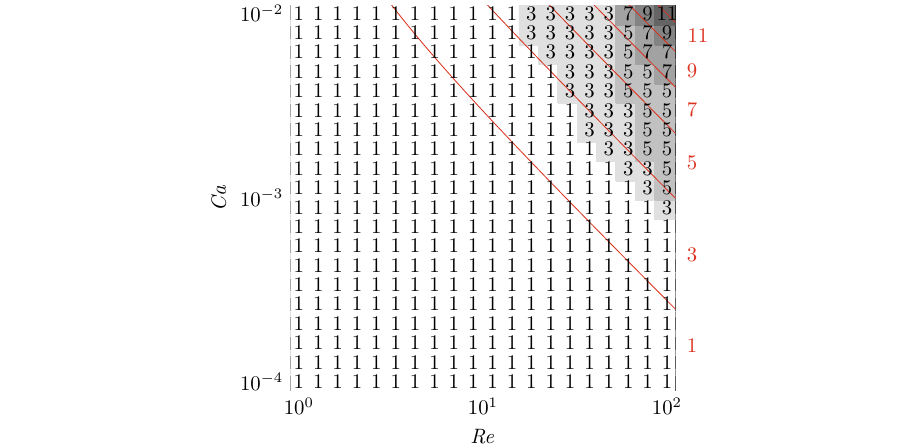}%
    \caption{The minimum number of actuators required to stabilise the Navier-Stokes film compared to the number of unstable modes of the linearised weighted-residual system (red). The number of controls needed to stabilise the uniform film never exceeds the number of unstable modes of the linear system \(n_u\) as given in \cref{eqn:unstable_modes}. The ranges for the two parameters cover a broad range of different fluids, select examples are listed in \cref{sec:parameter_values}. Videos of selected instances of film evolution and control are available as supplementary material.}%
    \label{fig:m_min}%
  \end{figure}

  In \cref{fig:m_min} we compare our predictions for \(n_u\) from expression~\cref{eqn:unstable_modes} to the minimum number of controls required to stabilise the film in our numerical experiments of the Navier-Stokes system as \(\Rey\) and \(\Ca\) vary. We see that, as expected, the system is stabilisable at \(M \geq n_u\) in all cases. In fact, in the majority of the parameter space, the minimum number of actuators required to stabilise the uniform state is lower than the number predicted by the linear analysis, particularly at lower Reynolds numbers. As previous work by Salamon, Armstrong, and Brown~\cite{salamon1994traveling} and Ruyer-Quil and Manneville~\cite{ruyer2000improved} shows that important physical characteristics such as travelling wave speed begin to diverge from DNS results at \(\Rey \approx 5\), the fact that controls based on a linearisation of these equations match (or even exceed) the expected performance up to \(\Rey\approx 100\) is remarkable. However, after this point it becomes clear that we are reaching the limit of the model's validity, and the ability of the controls to stabilise the uniform state becomes less predictable. We note that at larger Reynolds and capillary numbers the film takes much longer to respond to the effects of the controls, making it more challenging to assess whether the uniform state is stabilisable. By dynamically estimating the sign of the fitted damping rate we can avoid running simulations over unfeasibly long times.

\section{Conclusion}
  \label{sec:conclusion}
  The research presented herein has demonstrated new and significant capabilities in terms of design and analysis of optimal feedback controls for complex physical systems. The stabilisation of the canonical multi-scale framework of a thin liquid film falling down an inclined plane by employing reduced-order models such as the Benney and first-order weighted-residual equations has been used as the physical setup for our proposed methodology. We developed an LQR approach via blowing and suction controls which has been shown to outperform the predictions of linear stability theory, and can successfully function beyond the region of model validity for either the Benney- or the weighted-residual-derived controls. 

  We have shown that even the crude controls used here far exceed their expected performance, and this opens up numerous avenues for future work. It remains to be seen whether higher-order models such as the second-order weighted-residual integral boundary layer model proposed by Ruyer-Quil and Manneville~\cite{ruyer2000improved} can be used to further improve the type of control demonstrated here. In addition, it would be desirable to remove the control dependence on discretisation by developing infinite-dimensional controls, which might also allow for an improved analysis of control performance.

  Although here we have performed numerical experiments to showcase the control efficacy, physical experiments on real fluids are an obvious next step that we hope our work will inspire. In order to achieve this in practice there are a number of useful assumptions that must be relaxed, namely the 2D nature of the flow and periodic boundary condition formulation. The additional dimension will allow for cross-flow instabilities (an interaction which needs to be further quantified), and the boundaries can also affect the stability of the film~\cite{pollak2011side}. The blowing and suction controls used in the present work offer a valuable theoretical foundation permitting a comprehensive examination of control performance for this system. We envision realistic embodiments thereof to require further analysis. Nevertheless, the developed methodological platform offers a promising springboard for both mathematical progress and transfer towards other forms of actuation within related control mechanisms. Finally, we recognise that the assumption that full observations of the interfacial height are available is often unrealistic. In these scenarios, adaptations of the LQR method such as static and dynamic output feedback controls have been used to stabilise long-wave models~\cite{thompson2016stabilising}, and so we are hopeful that future methods underpinned by the present work will generalise to the full Navier-Stokes system, and further to physical experiments.

\appendix
\section{Parameter values}
  \label{sec:parameter_values}
  Although the majority of the results in this paper are applied to the dimensionless systems governed by the dimensionless numbers \(L\), \(\theta\), \(\Rey\), and \(\Ca\), it is important not to forget the physical roots of these systems. For all of the numerical simulations in this work, we have fixed the aspect ratio \(L = 30\), the inclination angle \(\theta = \pi/3\), gravitational acceleration \(g=9.807\unit{m}\unit{s}^{-2}\), and control width \(\omega = 0.1\). A range of values for the dimensional parameters (and the resulting dimensionless numbers) is provided in \cref{tab:params}. A wide range of physical configurations of interest are thus described by a parametric envelope given by \(10^0 < \Rey < 10^2\) and \(10^{-4} < \Ca < 10^{-2}\).

  \begin{table*}[ht]
    \centering%
    \begin{tabular}{l|ccc|cc}
      \hline
      Fluid & \(\rho\) (\(\unit{kg}\unit{m}^{-3}\)) & \(\mu\) (\(\unit{kg}\unit{m}^{-1}\unit{s}^{-1}\)) & \(\gamma\) (\(\unit{N}\unit{m}^{-1}\)) & \(\Rey\) & \(\Ca\) \\ [0.5ex]
      \hline
      Water    & \(999.8\) & \(8.91\times10^{-4}\) & \(0.072\)  & \(28.2\) & \(0.0018\) \\
      Ethanol  & \(789.5\) & \(1.06\times10^{-3}\) & \(0.022\)  & \(12.6\) & \(0.0047\) \\
      Pentane  & \(626.0\) & \(2.24\times10^{-4}\) & \(0.018\)  & \(178\)  & \(0.0045\) \\
      Nitrogen & \(3.44\)  & \(6.88\times10^{-6}\) & \(0.0085\) & \(5.69\) & \(5.26\times10^{-5}\) \\
      \hline
    \end{tabular}%
    \caption{Parameters (and resulting dimensionless numbers) for a range of physical fluids with a Nusselt film height of \(175\times10^{-6}\unit{m}\).}%
    \label{tab:params}%
  \end{table*}


\section{Numerical simulations}
  \label{sec:numerical_simulations}
  The Navier-Stokes equations (\cref{eqn:momentum1,eqn:momentum2,eqn:wall,eqn:stress1,eqn:stress2}) are solved on a finite domain \(\Omega = [0,L]\times[0,8]\) (the permissive height setup has been designed to prevent spurious pressure waves in the gas affecting the film) using the volume-of-fluid (VOF) method~\cite{scardovelli1999}. The computations were performed using Basilisk~\cite{basilisk}, a free extension to the C language designed to simplify writing code to numerically solve PDEs. It solves the incompressible Navier-Stokes equations on an adaptive quadtree grid~\cite{popinet2015} using the Bell-Collela-Glaz advection scheme with a CFL-limited time step, and an implicit viscosity solver (as did its predecessor, Gerris~\cite{popinet2003,popinet2009}). The grid spacing ranges from \(L\times 2^{-8}\) (covering the liquid film) to \(L\times 2^{-6}\) (smoothing out spurious pressure waves in the gas at the top of the finite computational domain). The time step is restricted to a maximum of \(0.05\) to prevent sudden jumps in the actuator inputs.

  Since the control strategy is fundamentally agnostic to the specifics of the PDE system being controlled aside from the entries of the linearised matrices \(J\) and \(\Psi\), the control code can be largely separated from the fluid simulation code. It would thus be relatively easy to transfer the same framework to a different problem.

  The Benney and weighted-residual equations are solved using second-order finite-difference stencils for the spacial grid and a second-order backward finite-difference scheme (BDF2) in time as, in Thompson, Tseluiko, and Papageorgiou~\cite{thompson2016falling}. The resulting problem is fully implicit and is solved via direct Newton iteration. All the computations in this paper were performed on a grid with a spacing of \(L\times 2^{-8}\) to match the resolution of the Basilisk grid.

\section{Control theory fundamentals}
  \label{sec:control_theory_fundamentals}
  Here we provide a brief overview of some important definitions in control theory relevant to our study. For more detailed aspects we refer the interested reader to the seminal work of Zabczyk~\cite{zabczyk2020mathematical}.

  Suppose we have the linear control system
  \begin{equation}
    \label{eqn:cont_system}\hat{h}_t = J\hat{h} + \Psi u, \qquad u = K y, \qquad y = \Phi\hat{h},
  \end{equation}
  which can be written \(\hat{h}_t = (J + \Psi K\Phi)\hat{h}\). The pair \((J,\Psi)\) is \textit{controllable} if, for any pair of states \(\hat{h}_0,\hat{h}_1\in\mathbb{R}^N\) there exists a control \(u\) that takes \(\hat{h}\) from \(\hat{h}_0\) to \(\hat{h}_1\) in finite time. The pair \((J,\Phi)\) is \textit{observable} if for all initial conditions \(\hat{h}_0\in\mathbb{R}^N\) there exists a time \(T>0\) after which \(\hat{h}_0\) is uniquely determined from the observations \(\{y(t) \vert t\in[0,T]\}\). Controllability and observability are duals, that is, if \((J,\Psi)\) is controllable then \((J^*,\Psi^*)\) (where \(\cdot^*\) is the conjugate transpose) is observable and conversely if \((J,\Phi)\) is observable then \((J^*,\Phi^*)\) is controllable.

  We can check if a pair \((J,\Psi)\) is controllable with the \textit{Kalman rank condition}: \((J,\Psi)\) is controllable if \(\text{rank}([J\vert \Psi]) = N\), where
  \begin{equation}
    [J\vert \Psi ] = [\Psi \ J\Psi \ J^2\Psi \ \ldots\ J^{N-1}\Psi ]
  \end{equation}
  is known as the controllability matrix.

  In this paper, we are concerned with controlling towards the state \(\hat{h} = 0\) rather than an arbitrary interface (see Thompson et al.~\cite{thompson2016stabilising}), and so we require a weaker form of controllability. For this we require \((J,\Psi)\) to be \textit{stabilisable}, which means that there exists a gain matrix \(K\) such that \(J+\Psi K\) is stable (i.e.\ has strictly negative real parts to all its eigenvalues). Similarly, in this case \((J,\Phi)\) is \textit{detectable} if we can choose an \(L\) such that \(J + L\Phi\) is stable, corresponding to being able to observe all of the unstable modes of the system. As for controllability and observability, stabilisability and detectability are dual properties (simply set \(L=K^*\) and vice versa).

\rnew{
\section{Initial Conditions}
  \label{sec:initial_conditions}
  It is not immediately obvious that a single numerical experiment starting from a travelling wave solution of the Navier-Stokes equations is sufficient to assess controllability at a given point in the parameter space. However, it turns out that this single initial condition captures the most challenging dynamics of the system. This is because, for the range of parameters used here, it is known that the observed long time behaviour is a travelling wave solution corresponding to the most unstable mode, which was observed both numerically and experimentally by a wide range of studies \cite{chang1994wave, denner2018solitary, kalliadasis2011falling, liu1993onset, nosoko2004evolution}. By following the uncontrolled evolution of a variety of randomly constructed initial conditions for the interface shape $h(x,0)$ in \cref{fig:random_ic} (left panel), we have numerically confirmed this wave to be independent of the initial condition, as showcased in \cref{fig:random_ic} (centre panel).

  \begin{figure}[ht]
    \centering%
    \includegraphics[width=\textwidth,page=1,trim={0 0 0 0},clip]{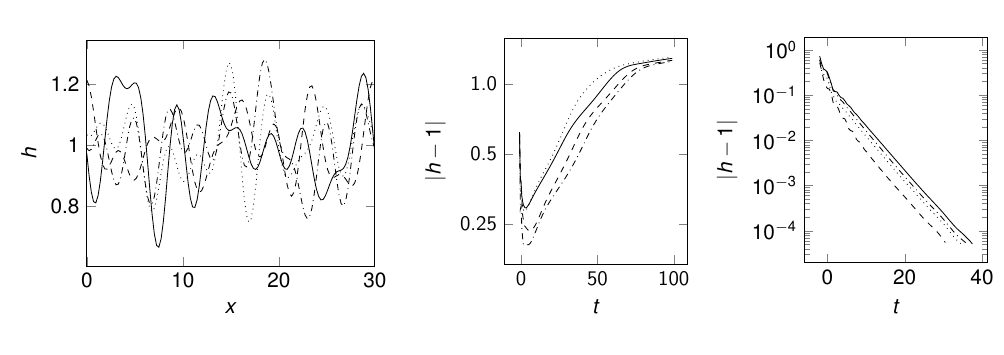}%
    \caption{\rnew{(left) Interface shape random initial conditions \(h(x, 0) = 1 + 0.1\sum_{m=1}^9 a_m \cos(mx + b_m)\), where \(a_m\) and \(b_m\) are chosen uniformly from \([0,1]\). (centre) If left uncontrolled, the stable modes quickly decay and the most unstable mode \(m=1\) dominates the dynamics, growing until nonlinear effects become significant and a travelling wave forms. (right) Controls switched on before the most unstable mode has had a chance to grow leading to a damping of perturbations with very similar rates.}}%
    \label{fig:random_ic}%
  \end{figure}

  In \cref{fig:random_ic} we also illustrate how the controls lead to very similar damping rates of the perturbation towards the target flat solution (right panel). This is intuitive because the most unstable mode has not had the chance to develop into its equilibrium state -- rather than needing to control a large amplitude wave (the fully developed travelling wave solution), for a random initial condition we have a set of perturbations with smaller magnitudes, some of which are stable even when uncontrolled.  For this reason, for all of the \textit{in silico} experiments applying the LQR controls to the numerical simulations of the full nonlinear Navier-Stokes simulations, we begin from the steady travelling wave that represents the most stringent test for our control strategy.
  
}

\section*{Supplementary material}
  \label{sec:supplementary}
  Supplementary material showing the evolution of the interface before and after the application of controls alongside the corresponding 2-norm deviations across a range of Reynolds and capillary numbers will be available upon publication.

  The version of the code used for this paper, along with installation instructions and documentation, can be found on \href{https://github.com/OaHolroyd/falling-film-control/tree/paper-dec-2022}{GitHub}. On a single core a full simulation (for instance the one shown in \cref{fig:controlled_interfaces}) takes \(\sim 10\) hours for the Navier-Stokes and weighted-residual systems (the Benney system is considerably faster).

\section*{Acknowledgements}
  \label{sec:acknowledgements}
  Oscar Holroyd is grateful for the computing resources supplied by the University of Warwick Scientific Computing Research Technology Platform (SCRTP) and funding from the UK Engineering and Physical Sciences Research Council (EPSRC) grant EP/S022848/1 for the University of Warwick Centre for Doctoral Training in Modelling of Heterogeneous Systems (HetSys). Radu Cimpeanu and Susana Gomes also acknowledge EPSRC support via grant EP/V051385/1. For the purpose of open access, the authors have applied a Creative Commons Attribution (CC BY) licence to any arising Author Accepted Manuscript version.

\bibliographystyle{siamplain}
\bibliography{references}
\end{document}